\documentclass{amsart}

\usepackage{amssymb}
\usepackage[all]{xy}
\usepackage{hyperref}

\usepackage{enumitem}   

\setlist[enumerate]{itemsep=.2em,topsep=.2em,leftmargin=1.25em,itemindent=2.0em}

\newtheorem{thm}{Theorem}

\newtheorem{lem}[thm]{Lemma}
\newtheorem{cor}[thm]{Corollary}

\newtheorem{prop}[thm]{Proposition}

\newtheorem{conj}[thm]{Conjecture}

\theoremstyle{definition}
\newtheorem{defn}[thm]{Definition}

\newtheorem{say}[thm]{}
\newtheorem{exmp}[thm]{Example}

\newtheorem*{ack}{Acknowledgments}      
\newtheorem{notation}[thm]{Notation}   
  
\newtheorem{defn-thm}[thm]{Definition--Theorem}  
\newtheorem{defn-lem}[thm]{Definition--Lemma}  

\theoremstyle{remark}

\renewcommand{\c}[0]{{\mathbb C}}

\newcommand{\z}[0]{{\mathbb Z}}

\renewcommand{\r}[0]{{\mathbb R}} 

\renewcommand{\a}[0]{{\mathbb A}}

\newcommand{\p}[0]{{\mathbb P}}
\newcommand{\f}[0]{{\mathbb F}}
\newcommand{\q}[0]{{\mathbb Q}}
\newcommand{\map}[0]{\dasharrow}
\newcommand{\qtq}[1]{\quad\mbox{#1}\quad}

\newcommand{\mult}[0]{\operatorname{mult}}
\newcommand{\discrep}[0]{\operatorname{discrep}}

\newcommand{\aut}[0]{\operatorname{Aut}}

\newcommand{\sing}[0]{\operatorname{Sing}}

\newcommand{\cl}[0]{\operatorname{Cl}}

\newcommand{\simq}[0]{\sim_{\q}}

\newcommand{\bir}[0]{\operatorname{Bir}}

\def\into{\DOTSB\lhook\joinrel\to}

\def\loccoh#1.#2.#3.#4.{H^{#1}_{#2}(#3,#4)}

\DeclareMathAlphabet{\mathchanc}{OT1}{pzc}%
                                {m}{it}

\newcommand{\gm}[0]{{\mathbb G}_m}

\newcommand{\GL}{\mathrm{GL}}
\newcommand{\PGL}{\mathrm{PGL}}

\usepackage[all]{xy}\xyoption{dvips}

\begin{document}
\bibliographystyle{amsalpha}


 \title[Cubic surfaces]{Cubic surfaces with infinite,\\ discrete automorphism group}
 \author{J\'anos Koll\'ar and David Villalobos-Paz}

       \begin{abstract}  We prove that the automorphism group of an affine, cubic surface with equation     $xyz=g(x,y)$ contains  ${\mathbb Z}$ as a finite index subgroup. These equations  were first studied by Jacobsthal (1939) and Mordell (1952).
       \end{abstract}

              \maketitle

While the automorphism group of 
       a smooth,  projective, cubic surface $\bar S$ is always finite, 
       the automorphism group of 
       a smooth,  affine, cubic surface $S$ can be  infinite.  The Markov-type equations
       $xyz=x^2+y^2+z^2+c$  give the best known examples.

       In general,   non-linear automorphisms of  a  smooth, cubic surface  $S\subset \a^3$ seem to be controlled by the singularities of the pair $(\bar S, \bar B )$, where $ \bar B :=\bar S\setminus S$ is the {\it curve at infinity.}
              If $ \bar B $ is smooth
        then every  automorphism  of $S$ is linear, and 
        $\aut(S)$ is finite.
        If $ \bar B $ has a unique singular point and $\bar S$ is smooth,
        then $\aut(S)$ is finite, but contains a non-linear involution; see    \cite[18]{k-lk3}.
However, if $ \bar B $ has 2 or 3 singular points and $\bar S$ is smooth,
then $\aut(S)$ is infinite.  These are the  Markov-type surfaces, see  (\ref{mar.typ.defn}).
No other smooth, affine,  cubic surface with infinite, discrete automorphism group seems to have been known.  


For cubic surfaces $S_g:=\bigl(xyz=g(x,y)\bigr)$,  the  curve at infinity is an irreducible, nodal cubic, and $\bar S_g$ has an $A_2$ singularity at the node. 
We show that, over any field,  all of these have an infinite automorphism group, and describe all
isomorphisms between them. The answer---given in  (\ref{main.isom.thm})---turns out to be quite subtle,  and less complete  over $\z$.  The following are representative examples.

\begin{exmp} \label{S0.main.thm}  For $S_0:=\bigl(xyz=x^3+y^3+1\bigr)\subset \a^3_{\z}$ the following hold.
     \smallskip

    {\rm (\ref{S0.main.thm}.1 Generators.)}   $\aut_{\z}(S_0)$ is generated by 3 involutions
    $$
    \begin{array}{rl}
      \tau\ : &(x,y,z)\mapsto (y,x,z),\\
  \sigma_x:& (x,y,z)\mapsto \bigl(x, xz-y^2, xz^2-y^2z-x^2y\bigr),\qtq{and}\\
  \sigma_y:&(x,y,z)\mapsto \bigl(yz-x^2, y, yz^2-x^2z-xy^2\bigr).
    \end{array}
    $$
    
    {\rm (\ref{S0.main.thm}.2 Structure.)}   $\aut_\z(S_0)\cong \z\rtimes \z/2$.

        \smallskip
 {\rm (\ref{S0.main.thm}.3 Fundamental domain.)} 
  $\{(x,y,z): x\leq y\leq  (x^3+1)^{1/2}\}$ is a fundamental domain for the action of $\aut_{\z}(S_0)$ on  $S_0\cap (\r^+)^3$.
\end{exmp}

\begin{exmp} \label{S12.main.thm}  Set 
$T_{n,m}:=(xyz=x^3+y^3+nx+my+1)$ for $n, m\in \z\setminus\{0\}$. Then 
  $$
  \aut_\z(T_{n,m})=\aut_\c(T_{n,m})\cong
  \begin{cases}
    \z\qtq{if} n\neq m, \qtq{and}\\
    \z\rtimes \z/2 \qtq{if} n= m.
    \end{cases}
  $$
  In their simplest form, 
the  coordinate functions  of the generator of $\z$ are given by
polynomials of degrees $13$,  $34$, and  $55$,   containing $110$, $998$, and  $2881$   monomials.
\end{exmp}

We also study  equations of the form
$xyz=g(x,y)$ for higher degree polynomials $g$, which were  
 considered  in  \cite{MR1557037} and \cite{MR51852}.

\begin{thm} \label{Sg.main.thm}
  Let $R$ be an integral domain, $g(x,y)\in R[x,y]$, and set
  $$
S_g:=\bigl(xyz=g(x,y)\bigr)\subset \a^3_R.
\eqno{(\ref{Sg.main.thm}.1)}
$$
Assume that $g$ has bidegree $(n, m)$ with $n, m\geq 3$, and the coefficients of $x^n, y^m, 1$ are units in $R$. 
Then there is an explicit,  infinite order automorphism $\sigma_g$ of $S_g$, defined over $R$, to be constructed  in (\ref{gen.not.alg.say}--\ref{Sg.main.thm.pf}).
\end{thm}

\begin{say}[Properties of $\sigma_g$]
  If $g$ has bidegree $(3,3)$ then no 1-dimensional, closed subscheme of $S_g$ is invariant under $\sigma_g$. Thus  $\sigma_g$ is loxodromic (as in \cite{cantat2024}) and in many cases this implies that integral points are Zariski dense; see (\ref{fin.orb.say}). Most likely the same holds for all $n,m\geq 3$.

  By contrast, if $g$ has bidegree $(2,2)$, then
  $S_g$ is a Markov type surface  (\ref{mar.typ.defn}) and  $\sigma_g$ preserves the pencil of curves $(\lambda x+\mu y=0)$, thus  $\sigma_g$ is parabolic.
  Integral points are usually not  Zariski dense; see
\cite[p.297]{MR249355}.

  We believe that $\langle \sigma_g\rangle$ has finite index in $\aut_R(S_g)$, and in most cases $\aut_R(S_g)=\langle \sigma_g\rangle$. We
prove these for  cubic polynomials in  (\ref{gen.aut.cor.exmp}) and  (\ref{gen.aut.cor.2}).
  \end{say}

\begin{say}[Steps of the proof]
For Example~\ref{S0.main.thm} the involutions are constructed in
(\ref{comp.sig.S0.say}). It is then easy to see that the subgroup
$\langle \sigma_x, \sigma_y, \tau\rangle$ is 
isomorphic to  $D_{\infty}$  (\ref{struct.sg.say}), and its
fundamental domain is as claimed (\ref{fun.dom.say}).
Most of the paper is then devoted to proving that there are no other automorphisms. This is obtained in (\ref{gen.aut.cor.3}), as  a special case  of (\ref{main.isom.thm}), which studies the surfaces (\ref{Sg.main.thm}.1) for
cubic polynomials $g(x,y)$.
The proof of (\ref{main.isom.thm}) is outlined in
(\ref{main.isom.thm.plan}), and completed in 
(\ref{main.isom.thm.pf}).

For 
Example~\ref{S12.main.thm} and Theorem~\ref{Sg.main.thm},
we construct the analogs of
$\sigma_x, \sigma_y$ in 
(\ref{gen.not.alg.say}),  and prove Theorem~\ref{Sg.main.thm} in (\ref{Sg.main.thm.pf}).
We have a complete description of all isomorphisms only for $\deg g=3$, see  (\ref{main.isom.thm}).

\end{say}
   
        \begin{say}[Previous results on automorphisms]\label{prev.res.say}
       \cite[4.2.2]{MR0769779}   showed the existence of   infinite order automorphisms of
       $P_0:=\p^2\setminus(xyz=x^3+y^3)$.
              Explicit, degree 8 polynomials  defining these automorphisms
              are given in  \cite[Sec.6]{MR1942244}.
               These maps  were 
               generalized to other Del~Pezzo surfaces in \cite{mcduff2024singular}.

       It was  observed in  \cite{k-lk3}  that  suitable  variants of these maps are Geiser-type involutions (\ref{g.inv.pencil.say}), occurring in pairs  $\sigma_{\pm}$.
       Since $\sigma_{\pm}$ are involutions, 
       $$
       \sigma_{\pm}^2(x{:}y{:}z)= \Phi_{\pm}(x,y,z)\cdot (x{:}y{:}z),
              $$
       for some polynomials $\Phi_{\pm}$. 
       Working over $\z$ suggests that the $ \Phi_{\pm}$ should be powers of $xyz-x^3-y^3$. 
       We computed, and indeed 
       $$
       \Phi_{\pm}(x,y,z)=(xyz-x^3-y^3)^{21}. 
       $$
       The  integral points of  $P_0$
             are  given by
       $(x,y,z)\in \z^3$ that satisfy $xyz-x^3-y^3=\pm 1$.
       Thus the $\sigma_{\pm}$ map integral points of $P_0$
       to integral points.

       The cubic surface $S_0$ in Example~\ref{S0.main.thm} is the universal cover of $P_0$. The covering map, given by
       $(x,y,z)\mapsto (x{:}y{:}z)$, has degree 3, and 
       the   $\sigma_{\pm}$  lift to automorphisms of $S_0$.
       Surprisingly,  in the latter incarnation one can take a `twisted square root' of the   $\sigma_{\pm}$; these are the automorphisms  $\sigma_x,\sigma_y$
        in (\ref{S0.main.thm}.1). By explicit computation, 
       $$
       \sigma_+=\sigma_x\circ\tau\circ \sigma_x, \qtq{and}
       \sigma_-=\sigma_y\circ\tau\circ \sigma_y.
       \eqno{(\ref{prev.res.say}.1)}
       $$
       
       Building on these results,  \cite{k-lk3} classifies all 
       smooth Del~Pezzo surfaces $\bar S$ and nodal curves $\bar B\subset \bar S$
       such that $\aut(\bar S\setminus \bar B)$ is infinite. In all cases $\deg \bar S\geq 4$.

       In many papers, the question of unexpected isomorphisms is  approached from the other direction, as the search for unexpected compactifications.
       The best studied example  is  $\a^2$.
              For compactifications  with Du~Val singularities see
\cite{MR0961326}, while 
  \cite{sawahara2024}  classifies even more general cases.
          \end{say}

        The constructions can be generalized to give   isomorphisms between some other affine Del~Pezzo surfaces, and    exclude   other possibilities.
Currently not all cases are covered, but we consider the following likely.

   \begin{conj} Let $S$ be a  smooth, affine, cubic surface whose  automorphism group is  infinite and discrete. Then $S$  is
          \begin{enumerate}
          \item either a Markov-type surface (\ref{mar.typ.defn}),
          \item or the curve at infinity $\bar B$ is nodal, and $\bar S$ is singular at the node.
          \end{enumerate}
   \end{conj}

\subsection*{Integral points}{\ }

\cite{MR1557037} studied  $xyz=g(x,y)$ when $g$ is symmetric in $x,y$. 
\cite{MR51852} shows that  $xyz=ax^3+by^3+c$ always has  integer solutions, for example
$$
x=bu^8+3abc^2u^5+3a^2bc^4u^2+c, \quad
y=u^3+ac^2,
\qtq{with} u:=a^3bc^5+1.
$$
For $S_0$ this gives  $( 365, 9, 14803)$.

\cite{MR51852} and \cite[Sec.30.3]{MR249355} state that the same holds for all
$xyz=g(x,y)$. This turns out to be not completely true; see \cite{MR3320468, MR3826636} for corrections.



For the growth rate of the number of solutions, we propose the following.

\begin{conj}\label{int.pts.conj}
  The number of integral points on $S_0$ satisfying
  $1\leq x,y\leq N$ grows like $c\log^2 N$.
\end{conj}

Numerical searches suggest that the true value of $c$ should be close to 1.
We prove in (\ref{get.log2.gr.say}.5) that the  number of integral points   is at least  $\epsilon\log N\log\log N$  for some (very small) $\epsilon>0$.
 
For the upper bound, note first that the number of integral points  satisfying
$-N\leq x,y\leq N$ is asymptotically $6N$, but almost all of these points lie on 3 lines, see (\ref{aff.l.say}). If we remove these lines, the
number of  points should grow like $cN^{1/2}$, but again most of these points should lie on some other curves $\a^1\to S_0$.

Thus the logarithmic upper bound  in Conjecture~\ref{int.pts.conj} implies that,
for every $\a^1\to S_0$ defined over $\z$, the image has finite intersection with $ (\r^+)^3$. 
In all examples that we know, the image is  disjoint from   $ (\r^+)^3$;
see Section~\ref{int.pts.sec} for details.

For all equations as in (\ref{Sg.main.thm}.1), \cite{vanluijk, MR3732687, browning2024integralpointscubicsurfaces} suggest that   the number of integral points on $S_g$, not lying on the image of any curve  $\a^1\to S_g$, and satisfying
$-N\leq x,y\leq N$, should grow like  a power of $\log N$.

Unfortunately, Conjecture~\ref{int.pts.conj} and
\cite{vanluijk,  browning2024integralpointscubicsurfaces}
suggest 3 different exponents for  $\log N$; see
(\ref{which.power.rem}) for a discussion.


Composing the lines in $S_0$ with $\aut_\z(S_0)$  gives infinitely many $\a^1\to S_0$ defined over $\z$, given by higher and higher degree polynomials.
Thus we need to remove infinitely many curves in order to have a chance for 
logarithmic growth rate if negative values of $x,y$ are allowed.


        \begin{ack} 
          We thank    T.~Browning, P.~Hacking,  D.~McDuff, J.~Li,   K.~Oguiso, and  K.~Siegel  for  many  helpful comments;  M.~Bennett, B.~Grechuk,  P.~Sarnak  and C.~Skinner for   references to the number theory literature. 
          The computations were done using Macaulay2 \cite{M2}. 
  
  Partial  financial support  to JK   was provided  by  the NSF under grant number
DMS-1901855 and by the Simons Foundation   under grant number SFI-MPS-MOV-00006719-02.
        \end{ack}

        \section{Cubic surfaces}

\begin{notation}\label{s.b.noit.not}
From now on  $S:=\bigl(h(x,y, z)=0\bigr)\subset \a^3$ denotes an affine cubic surface. Its projective closure is denoted by 
$\bar S\subset \p^3$, and  $ \bar B :=\bar S\setminus S$ is the {\it curve at infinity.}   It is  defined by the homogeneous cubic part $h_3(x,y, z)$ of $h(x,y, z)$.

To avoid  degenerate cases, we assume that  $\bar S$ has at worst  Du~Val singularities. Equivalently, $\bar S$ is normal and not a cone.
\end{notation}

       \begin{say}[Linear automorphisms]\label{lin.aut.say}
         For an affine cubic surface $S$, the group of
         {\it linear automorphisms,}
         denoted by $\aut^{\rm lin}(S)$, is the largest subgroup of
         $\aut(\bar S)$ that maps $S$ to itself. Equivalently, maps the
         plane at infinity to itself. This is thus a subgroup of $\PGL_4$.

         The automorphism group of a smooth, projective cubic surface is finite; all possibilities are listed in 
         \cite[Sec.9.5]{MR2964027}.

         If $\bar S$ has Du~Val singularities, then  $\aut(\bar S)$ can be   positive dimensional. All these are listed in   \cite{MR4206439}.
There are 9 isolated cubics and a 1-dimensional family. For example,   the surface 
$S_{14}:=(  x_0x_2^2 = x_1^3 + x_3x_0^2)$ has  an $E_6$ singularity at $(0{:}0{:}0{:}1)$, and  admits a $\gm$-action
$$
(x_0{:}x_1{:}x_2{:}x_3)\mapsto (x_0{:}\lambda^2x_1{:}\lambda^3x_2{:}\lambda^6x_3).
$$
Thus all 3 of the affine surfaces  $S_{14}\setminus (x_i=0)$  (for $i=0,1,2$) are
smooth with a $\gm$-action. As another example,
$(xyz=1)$ is isomorphic to  $\gm^2$, and its automorphism group is $\gm^2\rtimes \GL_2(\z)$.

For $S_0:=(xyz=x^3+y^3+1)$
the obvious automorphism is the involution $\tau:(x:y:z)\mapsto (y:x:z)$.

If $k$ contains a 3rd root of unity  $\epsilon\neq 1$, then
$\aut^{\rm lin}(S_0)$  also has a subgroup  $\z/3\oplus \z/3$
generated by  the diagonal matrices
$\mu=(\epsilon, \epsilon, \epsilon)$ and $\rho:=(\epsilon, \epsilon^2, 1)$.
In this case $\aut^{\rm lin}(S_0)$ 
has order 18.
 \end{say}

       \begin{say}[Markov-type surfaces]\label{mar.typ.defn}   Let $S\subset \a^3$ be a  cubic surface with projective closure
         $\bar S\subset \p^3$. Let $p\in \bar S\setminus S$ be a point at infinity.  Projecting $S$ from $p$ gives a birational map $\pi_p:\bar S\map \p^2$ if $p$ is a double point, and  a double cover if $p$ is a smooth point.
         In the latter case let $\tau_p: \bar S\map \bar S$ be the corresponding Galois involution.  Note that $\tau_p$ contracts
         $\bar S\cap T_p\bar S$ to $p$, and  $\tau_p$ gives an
         automorphism of $S$ iff the tangent plane $T_p\bar S$ is
         the plane at infinity. Equivalently, iff $p$ is a singular point of
         $\bar B$.

         If $p=(0{:}0{:}1{:}0)$ then the equation of $S$ can be written as
         $z^2+c_2(x,y)z+c_3(x,y)=0$, where $c_i$ has degree $\leq i$.
         Then
         $$
         \tau_p: (x,y,z)\mapsto \bigl(x,y, c_2(x,y)-z\bigr).
         \eqno{(\ref{mar.typ.defn}.1)}
         $$
         This is a   nonlinear automorphism if
         $\deg c_2=2$, but a linear one if $\deg c_2\leq 1$.

We call $S$ a {\it Markov-type} surface if $\bar B$ has 2 or 3 nodes, and $\bar S$ is smooth at  2 or 3 of them. (Other authors have different definitions.)

         For the classical  Markov equation
         $3xyz=x^2+y^2+z^2$ discussed in \cite{mar-1, mar-2}, the curve at infinity is $(xyz=0)$, with 3 singular points. Each of them gives an involution; they generate  a subgroup isomorphic to $\z/2*\z/2*\z/2$. The automorphims group seems to have been fully determined only in \cite{MR0342518, MR4337481}.
         Automorphism groups and integral points of  Markov-type surfaces 
         are  further  studied in  \cite{MR2649343, MR4448994}.
         
                  If the  curve at infinity is the union of a line and a conic meeting at 2 points, then we get
         the free product of 2 involutions.   
         In all other cases the curve at infinity has at most 1 singular point, giving a single,  order 2 subgroup of $\aut(S)$.

          \end{say}

       \section{Cubic surfaces with an $A_2$ singularity}

       \begin{say} \label{our.cubics.defn} Let $R$ be an integral domain with quotient field $K$.
         We study cubic surfaces
         $S\subset \a^3_R$ such that
         \begin{enumerate}
         \item the curve at infinity $\bar B_K$ is geometrically irreducible,  nodal, and
         \item  $\bar S_K$ has an $A_2$ singularity at the node.
         \end{enumerate}
         In suitable coordinates  the equation of $\bar S$ is
         $g_2(x,y)z=g_3(x,y,w)$, where $g_i$ is homogeneous of degree $i$.
         Then (\ref{our.cubics.defn}.2) holds iff $g_3(0,0,1)\neq 0$
         and (\ref{our.cubics.defn}.1) holds $g_2(x,y)$ and $g_3(x,y,0)$
         are relatively prime.
         \end{say}

       \begin{notation} \label{a(t).not}
         Let $R$ be an integral domain.
          We let  $P_n(R)$ denote the set of polynomials
          $a(t):=a_nt^n+\cdots+a_0$  satisfying
          $a_n\neq 0\neq a_0$.
          Two such polynomials $a(t)$ and $b(t)$  are {\it $\gm$-equivalent} iff  $b(t)=\lambda_1a(\lambda_2 t)$ for  some
          $\lambda_1, \lambda_2\in R^\times$.

          Let  $P^\times_n(R)\subset  P_n(R)$ denote the subset where $a_n, a_0$ are units. Thus  $P^\times_n(R)= P_n(R)$ if $R$ is a field.

       \end{notation}

       \begin{defn} \label{split.defn}
         Let $k$ be a field and consider a singularity  $Z$ defined by an equation
         $g_2(x_1,\dots, x_n)+(\mbox{higher terms})$ where
         $g_2$ is a homogeneous quadric of rank 2. We say that $Z$ is {\it split} if $g_2$ is the product of 2 linear forms over $k$.
         We use this notion for nodes ($n=2$) and for  $A_m$ singularities
         ($n=3$ and $m\geq 2$).
         \end{defn}

       \begin{say}[Normal form for the equation]\label{norm.f.lem}
         Let $k$ be a field and 
         $\bar S\subset \p^3_{xyzw}$ a cubic surface as in  (\ref{our.cubics.defn}).
         
If the $A_2$ singularity is split, we can choose $g_2=xy$. 
Then, by changing $z$ we can eliminte from $g_3$ all terms divisible by $xy$. That is, we have a cubic form $g_3$ restricted to $(xy=0)$.
We write these as
        $$
        (a_3x^3+a_2x^2w+a_1xw^2+a_0w^3)+(b_3y^3+b_2y^2w+b_1yw^2+b_0w^3)-cw^3,
        $$
where  $c=a_0=b_0$.  As we noted in  (\ref{our.cubics.defn}),
$a_3, b_3, a_0, b_0$ are all nonzero.
If $S$ is smooth then $a(x)$ and $b(y)$  have no multiple roots.
        We have thus proved the following.

\medskip
    {\it Claim \ref{norm.f.lem}.1.}  Let $k$ be a field and $S$ an affine cubic surface as in  (\ref{our.cubics.defn}) with  a split  $A_2$ singularity at infinity. Then  one can  write its  equation in the form 
          $$
          S_{a,b}:=\bigl(xyz=a(x)+b(y)-c\bigr)\subset \a^3,
        $$
      where
        $a(t), b(t)\in P_3(k)$  and $a(0)=b(0)=c$. \qed
\medskip

      The only linear isomorphisms between these normal forms are given by 
      $$
      (x,y,z)\mapsto (\lambda_x x, \lambda_y y, \lambda_z z) \qtq{or}
      (\lambda_y y, \lambda_x x, \lambda_z z),
       $$
       for $\lambda_x, \lambda_y, \lambda_z\in k^{\times}$.
       \qed
       \end{say}

       Over other rings, we take  (\ref{norm.f.lem}.1) as our definition.

       \begin{defn} \label{Sab.def.R}
         Let $R$ be an integral domain, and $a(t),b(t)\in P^\times_3(R)$. Then $a_0=a(0)$ and $b_0=b(0)$. Set
         $$
          S_{a,b}:=\bigl(xyz=b_0a(x)+a_0b(y)-a_0b_0\bigr)\subset \a^3_R.
          \eqno{(\ref{Sab.def.R}.1)}
          $$
          Since $a_0, b_0$ are units, we can divide by $a_0b_0$, and  change  $z$ to $a_0b_0z$ to get the equivalent form
          $$
          S_{a,b}:=\bigl(xyz=a_0^{-1}a(x)+b_0^{-1}b(y)-1\bigr)\subset \a^3_R.
          \eqno{(\ref{Sab.def.R}.2)}
          $$
\end{defn}

        \begin{say}[Blow-up representation]\label{blw.up.say}
          Projecting from the singular point to the $(x{:}y{:}w)$-plane contracts the 6 lines
          $$
          \begin{array}{l}
          (x=b_3y^3+b_2y^2w+b_1yw^2+b_0w^3=0)\qtq{and}\\
            (y=a_3x^3+a_2x^2w+a_1xw^2+a_0w^3=0).
            \end{array}
          \eqno{(\ref{blw.up.say}.1)}
          $$
          Let $\beta_i$ be the roots of $b(t)$ and
          $\alpha_i$  the roots of $a(t)$.

          We can thus obtain  $S_{a,b}$ by starting with
          $\a^2_{uv}$, blow up the 6 points  $(\alpha_i, 0), (0,\beta_i)$,
          and then removing the birational transform of $(uv=0)$. The map is given by
          $$
          x=u^2v, \ y=uv^2,\ z=a(u)+b(v)-c,\ w=uv.
          $$
          If the 6 points  $\alpha_i, \beta_i$ are different, we get a smooth affine surface. If 2 (resp. 3) of them coincide, then the affine surface has an $A_1$  (resp.\ $A_2$) singularity. However, the singularity at infinity stays $A_2$.   These extra  $A_1$  or  $A_2$ singularities on $S_{a,b}$  play no role in our investigations.

          The  birational transform of $(w=0)$ becomes the curve at infinity $\bar B$.

          Thus the points   $\beta_{\infty}:=(0{:}1{:}0)$ and
          $\alpha_{\infty}:=(1{:}0{:}0)$ get identified to form the node at infinity.

          The 6 exceptional curves of the blow-up become the 6 lines through the singular point. 
The birational transforms of the lines $\langle \alpha_{\infty}, \beta_i\rangle$ or
$\langle \beta_{\infty}, \alpha_i\rangle$ become conics that are tangent to one of the branches of $\bar B$ at the node.
        \end{say}

        \begin{say}[Lines on $S_{a,b}$]\label{lines.Sab.say}
          There are 6 lines on $S_{a,b}$ passing through the singular point.
          If $\alpha$ is a root of $a(x)$, then $(\alpha, 0,t)$ is a line on 
          $S_{a,b}$. Similarly, if $\beta$ is a root of $b(x)$, then we get $(0, \beta,t)$.

          The residual intersection of  $S_{a,b}$ with  the plane spanned by these 2 lines is another line, giving 9  lines not  passing through the singular point.

          Thus we get a $k$-line on $S_{a,b}$   not  passing through the singular point iff both  $a(x)$ and $b(y)$ have a root in $k$.
          \end{say}

         \begin{say}[Moduli]\label{A2.mod.say} Using the blow-up representation (\ref{blw.up.say}), we see that the moduli of 
             projective cubics with an $A_2$ singularity is isomorphic to the
             moduli of $3+3$ points on a pair of lines  $(uv=0)\subset \p^2_{uvw}$.

             We are working with affine cubics; this corresponds to also fixing the line $(w=0)$. Thus we have  $(uv=0)\subset \p^2_{uvw}$ with marked points
             $(0{:}1{:}0)$ and  $(1{:}0{:}0)$. The automorphism group is now reduced to the maps
             $$
             (u{:}v{:}w)\mapsto (\lambda u{:}\mu v{:}w)\qtq{and}
             (\mu v{:}\lambda u{:}w).
             $$
             \end{say}

        \section{Constructing the involutions}

For  introductions   see \cite{geiser}, 
\cite[Chap.8]{MR2964027} or \cite[Sec.4]{k-lk3}.
        
        \begin{say}[Geiser involutions of singular Del~Pezzo surfaces]
          \label{g.inv.pencil.say}
          A degree 2 Del~Pezzo surface $T$ has a natural double cover structure
           $\pi:T\to \p^2$,  ramified along a degree 4 curve $R$. The 
          Geiser involution is the corresponding Galois involution.
          
          The preimage $E_L$ of a  line $L\subset \p^2$ is a  curve of arithmetic genus 1.
          
          Assume now that $T$ is singular at a point $p$.
Then $R$ is singular at $\pi(p)$. 
          Let $L\ni \pi(p)$ be a line. Then  $E_L$ is singular at $p$, we expect it to be  a rational curve with a node at $p$. The Geiser involution interchanges the 2 branches of $E_L$ at $p$. All such curves $E_L$ form a pencil in $|-K_T|$.
          This pencil is base point free on $T\setminus\{p\}$, its general fiber is $\gm$ and the special fibers, consisting of  pairs of lines, are
          isomorphic to $(uv=0)\subset \a^2_{uv}$.  The 
          Geiser involution interchanges the $u$ and $v$ axes.

          Our plan is to identify these pencils on $S_0$, and write down  formulas for $\sigma_x, \sigma_y$. Since $S_0$ is a  degree 3 Del~Pezzo surface, we need to blow it up at 1 point to get a  degree 2 Del~Pezzo surface.
          The only `natural' point on $\bar S_0$ is the singular point. So we look for pencils of rational curves whose base point is the singular point.
          Once we find these, it turns out that going from $\bar S_0$ to the  degree 2 Del~Pezzo surface involves the resolution of the singularity, 2 blow-ups, and 3 blow-downs.
          \end{say}
        
        \begin{say}[Computing $\sigma_x, \sigma_y$]\label{comp.sig.S0.say}
Projecting $S_0:=\bigl(  xyz=x^3+y^3+1\bigr)$  to the $(x,y)$-plane as in (\ref{blw.up.say})
shows that the open subset lying over   $(xy\neq 0)$  is isomorphic to   $\gm\times \gm$,
and the rest consists of 6 lines over the points
$(x=y^3+1=0)$ and $(y=x^3+1=0)$.

Thus projection to the  $x$ or $y$-axis gives an anticanonical pencil, whose
general fiber is $\gm$, and the special fibers consist of  pairs of lines.
The latter are over the points $(x^3+1=0)$ resp.\ $(y^3+1=0)$.

These projections are  good candidates for the pencils associated to the Geiser involutions in (\ref{g.inv.pencil.say}).
Thus, for the projection to the $y$-axis, in the $(x, y)$-coordinates we are looking for an involution of the form
$$
(x, y)\mapsto  \Bigl(\tfrac{c(y)}{x}, y\Bigr).
$$
For fixed $y$, this is an involution of $\gm\times \{y\}$ if $c(y)\neq 0$, but degenerates at the roots of $c(y)$. Since over the roots of $y^3+1$ we have reducible fibers,
it is reasonable to expect that $c(y)=y^3+1$.
This is a good guess, since $xyz=x^3+y^3+1$ gives that
$$
\tfrac{c(y)}{x}=\tfrac{y^3+1}{x}=yz-x^2.
$$
Thus we get a regular involution
             $$
                (x,y)\mapsto (\bar x, \bar y):=
\bigl(yz-x^2, y\bigr)=\Bigl(\tfrac{y^3+1}{x}, y\Bigr).
\eqno{(\ref{comp.sig.S0.say}.1)}
                $$
                We have no choice for $\bar z$, but it should be regular. 
               Using  that $\bar y^3+1=y^3+1=x\bar x$, we get that 
        $$
        \begin{array}{rcl}
        \bar z&=&\frac{\bar x^3+\bar y^3+1}{\bar x\bar y}=
        \frac{\bar x^3+x\bar x}{\bar x\bar y}=
        \frac{\bar x^2+x}{y}=\\[1ex]
        &=&
        \frac{y^2z^2-2x^2yz+x(x^3+1)}{y}=
        \frac{y^2z^2-2x^2yz+xy(xz-y^2)}{y}=
        yz^2-x^2z-xy^2.
        \end{array}
        \eqno{(\ref{comp.sig.S0.say}.2)}
        $$
        Thus on $S_0$ we obtain the   involutions
        $$
    \begin{array}{rl}
    \sigma_x:& (x,y,z)\mapsto \bigl(x, xz-y^2, xz^2-y^2z-x^2y\bigr)= \Bigl(x, \tfrac{x^3+1}{y}, \tfrac{(xz-y^2)^2+y}{x}\Bigr),\\[1ex]
  \sigma_y:&(x,y,z)\mapsto \bigl(yz-x^2, y, yz^2-x^2z-xy^2\bigr)=
   \Bigl(\tfrac{y^3+1}{x}, y, \tfrac{(yz-x^2)^2+x}{y}\Bigr).
    \end{array}
    \eqno{(\ref{comp.sig.S0.say}.3)}
    $$
 \end{say}

       \begin{say}[Structure of  $\aut_\z(S_0)$]\label{struct.sg.say}
         We determine the structure of the subgroup generated by $\sigma_x, \sigma_y$ and linear automorphisms.
         (We see in (\ref{main.isom.thm}) that this is in fact the whole  $\aut_\z(S_0)$.)
         We use the {\it infinite dihedral group} $D_{\infty}:=\z/2*\z/2$.

         First, $\sigma_x, \sigma_y$ are involutions, hence they generate a dihedral group. There are many ways to see that  $\sigma_x\sigma_y$ has infinite order; see for example (\ref{mult.red.prop}.1). So  $\langle \sigma_x, \sigma_y\rangle\cong D_{\infty}$.
         Next $\tau$ interchanges $\sigma_x, \sigma_y$, so
$$
   \bigl(\langle \sigma_x\rangle * \langle \sigma_y\rangle )\rtimes
  \langle \tau\rangle\cong \bigl(\z/2*\z/2)\rtimes \z/2.
  $$
  This is, however, also isomorphic to $D_{\infty}$, for example with generators
  $\sigma_x$ and $ \tau$.

  If $k$ contains a 3rd root of unity  $\epsilon\neq 1$, then, as noted in  (\ref{lin.aut.say}), 
$\aut^{\rm lin}(S_0)$  contains  a subgroup   $H$ of 
 order 9  generated by 
 $\mu=(\epsilon, \epsilon, \epsilon)$ and $\rho:=(\epsilon, \epsilon^2, 1)$. 
We have
$$
\sigma_x\mu=\rho\sigma_x, \ \sigma_x\rho=\mu\sigma_x, \ \tau\mu=\mu\tau, \
\tau\rho=\rho^2\tau.
$$
Thus $H$ is a normal subgroup.

Also note that 
the normalizer of  $\langle\mu\rangle$ descends to   $S_0/\langle\mu\rangle\cong P_0=\p^2\setminus(xyz=x^3+y^3)$. We see that 
$$
\sigma_x\tau\sigma_x\mu=\sigma_x\tau\rho\sigma_x=\sigma_x\rho^2\tau\sigma_x=
\mu^2\sigma_x\tau\sigma_x,
$$
so indeed $\sigma_x\tau\sigma_x$ descends. This was one of our starting points in  (\ref{prev.res.say}.1).
       \end{say}

        \begin{say}[Fundamental domain]\label{fun.dom.say}
          The rational forms of the involutions (\ref{comp.sig.S0.say}.3)
          $$
          \Bigl(x, \tfrac{x^3+1}{y}, \tfrac{(xz-y^2)^2+y}{x}\Bigr)\qtq{and}
          \Bigl( \tfrac{y^3+1}{x}, y, \tfrac{(yz-x^2)^2+x}{y}\Bigr)
          $$
           show that $\sigma_x, \sigma_y$ map  the positive octant $S_0\cap (\r^+)^3$ to itself. Applying $\tau$ if necessary, we may asume that $x\leq y$, and $(x^3+1)/y<y$ iff  $x^3+1< y^2$. Thus a fundamental domain  for the $\aut(S_0)$-action on  $S_0\cap (\r^+)^3$  is
              $$
              x\leq y\leq (x^3+1)^{1/2}.
              $$
        \end{say}

        Following the recipe of (\ref{comp.sig.S0.say}.2) for a general $S_{a,b}$ results in a $\bar z$ that is regular if $a_3=a_0$, but not otherwise.
        Next we discuss the geometric reason behind this. It will then show the correct generalization of the formulas (\ref{comp.sig.S0.say}.3).

        \begin{say}[Geometric construction of $\sigma_x, \sigma_y$]
          \label{geom.invol.say}
          First take the minimal resolution of the singularity of $\bar S_{a,b}$. We give 2 versions of the dual graph
          $$
          C_y {\mbox{ --- }} E_y  {\mbox{ --- }} E_x  {\mbox{ --- }}C_x
          \qtq{and}
          C_y {\mbox{ --- }} (-2)  {\mbox{ --- }} (-2)  {\mbox{ --- }}C_x
          $$
          where $E_x, E_y$ denote the exceptional curves, and
          $C_x, C_y$ the  birational transforms of the local branches of $\bar B$. The left side shows the names, the right side the self-intersections.
          
          Then we blow up twice the  point on $C_x$ that is also on the exceptional set. We get
        $$
\begin{array}{ccccc}
  E'_y & - & E'_x &  -  & E_1 \\
  \vert &&&& \vert\\
  C' && \rule[1pt]{40pt}{0.4pt} && {\ }E_2
\end{array}
\qtq{and}
\begin{array}{ccccc}
  (-2) & - &(-3) &  -  & (-2) \\
  \vert &&&& \vert\\
  (-1) && \rule[1pt]{40pt}{0.4pt} && {\ }(-1)
\end{array}
\eqno{(\ref{geom.invol.say}.1)}
$$
Thus (\ref{geom.invol.say}.1) has a numerical left-right symmetry.
The 2 blow-ups creating $E_1, E_2$ result in  a surface with  $(K^2)=1$, but it is not a weak Del~Pezzo surface since there is a $(-3)$-curve.

Note that the   $(-3)$-curve  $E'_x$ is the birational transform of $E_x$,
and by (\ref{blw.up.say}) it is intersected by the   birational transforms  of  3 lines, denote these by $L'_1, L'_2, L'_3$.

If we contract one of the lines $L'_i$, then 
(\ref{geom.invol.say}.1) becomes
$$
\begin{array}{ccccc}
  (-2) & - &(-2) &  -  & (-2) \\
  \vert &&&& \vert\\
  (-1) && \rule[1pt]{40pt}{0.4pt} && {\ }(-1)
\end{array}
\eqno{(\ref{geom.invol.say}.2)}
$$
Now we have a  degree 2 weak Del~Pezzo surface, and the
Geiser-type involution interchanges the left and right sides.

Remarkably, the end result is independent of which line we choose, but for now choose say $L'_3$.

To compute the  action of $\sigma_x$ on the birational transform $E''_x$ of $E_x$, we use the  coordinate $t=x/w$. Then  $E_y\cap E_x=(t=0)$ and
$E_x\cap C_x=(t=\infty)$.  The 3 lines  $L'_i$ meet $E_x$ at the 3 roots of 
$a(x)=a_3x^3+a_2x^2+a_1x+a_0$, call these $\alpha_i$. 

The  involution $\sigma_x$ then acts by
$t\mapsto c/t$ for some constant $c$.
The  Geiser involution maps lines to lines, thus it interchanges  $L''_1$ and  $ L''_2$, hence  it interchanges $\alpha_1$ and $ \alpha_2$.
This gives that  $c=\alpha_1\alpha_2$.  Thus
$$
\alpha_3\mapsto \tfrac{\alpha_1\alpha_2}{\alpha_3}.
$$
Note that $\sigma_x(\alpha_3)=\alpha_3$ iff $\alpha_1\alpha_2\alpha_3=1$.
That is, iff $a_3=a_0$.
This explains our stated observation that  the formula in  (\ref{comp.sig.S0.say}.2)
defines an involution if $a_3=a_0$, but not otherwise.

Once we understand that $\sigma_x$ need not be an automorphism of $S_{a,b}$, but an isomorphism between two different surfaces  $S_{a,b}$ and $S_{a',b'}$,
we choose the coordinates $t,t'$ on the surfaces independent of each other.
By (\ref{A2.mod.say})  they are well defined only up to scaling, thus we can as well choose
$t'=1/t$. 
That is, on the roots we get the map
$$
 (\alpha_1,\alpha_2,\alpha_3)\mapsto
\Bigl(\tfrac1{\alpha_1}, \tfrac1{\alpha_2}, \tfrac1{\alpha_3}\Bigr).
  $$
The inverses of the roots satisfy the equation
$\bar a(t)=0$,
where the coefficients  are in reversed order  (\ref{a(t).not}).
 \end{say}
  
        Building on these considerations---and after some trial  and error---the algebraic construction of the general
        $\sigma_x$ and $\sigma_y$ works out very cleanly.

 \begin{notation} \label{a(t).not.2}
         For  an integral domain $R$, let 
         $P_n(R)$ and $P^\times_n(R)$ be  as in (\ref{a(t).not}).
         For $a(t)\in P_n(R)$ define $a^*(t)$ by the formula
         $a(t)=ta^*(t)+a_0$. We also set
         $$
\bar a(t):=\tfrac1{a_0^{n-1}a_n}t^na\bigl(\tfrac{a_0}{t}\bigr),
$$
though the reason for this becomes clear only in (\ref{gen.not.alg.say}.8).

Note that  $a\mapsto \bar a$ is an involution on  $P^\times_n(R)$.
Furthermore, $\bar a_n=\tfrac1{a_0^{n-2}a_n}$ is in $R$ for $n\geq 3$
iff $a_n, a_0$ are units in $R$.

\end{notation}

\begin{say}[Algebraic construction of $\sigma_x$ and $\sigma_y$]\label{gen.not.alg.say}
  Start with a surface
  $$
  S_{a,b}=\bigl(xyz=a(x)+b(y)-c\bigr)\subset \a^3_{xyz},
  \eqno{(\ref{gen.not.alg.say}.1)}
  $$
  where $a(x)=a_nx^n+\cdots+a_0$ and $ b(y)=b_my^m+\cdots+b_0$. We assume that  $a_0=b_0=c$, and $a_n, b_m, c$ all nonzero.
As in (\ref{geom.invol.say})  
we want $\sigma_y$ to send $(x,y)$ to
$$
(\bar x, \bar y)=\bigl(b(y)/x, y\bigr)=(yz-a^*(x), y),
\eqno{(\ref{gen.not.alg.say}.2)}
$$
and we hope that the new surface has equation
$$
S_{\bar a,b}=\bigl(xyz=\bar a(x)+b(y)-c\bigr)
\eqno{(\ref{gen.not.alg.say}.3)}
$$
for some not yet known $\bar a(x)$, that must satisfy $\bar a_0=a_0$.

 Using that $b(y)=\bar x x$, we get that 
$$
\bar z=\frac{\bar a(\bar x)+b(\bar y)-c}{\bar x\bar y}=
\frac{\bar a(\bar x)-c+b(y)}{\bar x y}=
\frac{\bar x\bar a^*(\bar x)+x\bar x}{\bar x y}=
\frac{\bar a^*(\bar x)+x}{ y},
\eqno{(\ref{gen.not.alg.say}.4)}
$$
and it should be regular. 
Now use  that $\bar x=yz-a^*(x)$. So
$$
\bar a^*(\bar x)+x=\bar a^*\bigl(yz-a^*(x)\bigr)+x=
    (\mbox{divisible by $y$})+ \bar a^*\bigl(-a^*(x)\bigr)+x.
\eqno{(\ref{gen.not.alg.say}.5)}
$$
The only known polynomial in $x$ that is divisible by $y$ is
$$
a(x)=xyz-b(y)-c=y\bigl(xz-b^*(y)\bigr).
\eqno{(\ref{gen.not.alg.say}.6)}
$$
Thus our hope is that $a(x)$ divides  $\bar a^*\bigl(-a^*(x)\bigr)+x$.
A general $a(x)$ has distinct roots, so we need that 
 each root of $a(x)$ is a root of  $\bar a^*\bigl(-a^*(x)\bigr)+x$.  

 Let $\alpha$ be a root of $a(x)$.  Then
  $0=a(\alpha)=\alpha a^*(\alpha)+a_0$, so  $-a^*(\alpha)=a_0/\alpha$.
Substituting, we get that 
$$
\bar a^*\bigl(-a^*(\alpha)\bigr)+\alpha=\bar a^*\bigl(\tfrac{a_0}{\alpha}\bigr)+\alpha=
\tfrac{\alpha}{a_0}\bigl(\tfrac{a_0}{\alpha}\bar a^*\bigl(\tfrac{a_0}{\alpha}\bigr)+a_0\bigr)=\tfrac{\alpha}{a_0}\cdot \bar a\bigl(\tfrac{a_0}{\alpha}\bigr),
\eqno{(\ref{gen.not.alg.say}.7)}
$$
where we used that  $\bar a_0=a_0$.
Thus we are done if $\tfrac{a_0}{\alpha}$ is a root of $\bar a$.
This tells us that
$\bar a(x)$ should be a constant times $x^na\bigl(\tfrac{a_0}{x}\bigr)$.
Since $\bar a_0=a_0$, we must have
$$
\bar a(x):=\tfrac1{a_0^{n-1}a_n}x^na\bigl(\tfrac{a_0}{x}\bigr).
\eqno{(\ref{gen.not.alg.say}.8)}
$$
Thus, if $a_n, a_0$ are units, then $\sigma_y$ gives an $R$-isomorphism
between $S_{a,b}$ and $S_{\bar a, b}$.

If the $a_n, a_0$ are not units, then the
formulas (\ref{gen.not.alg.say}.2) for $(\bar x, \bar y)$ are in $R[x,y,z]$, but the formula for $\bar z$ involves the denominator $a_0^{n-2}a_n$.

Interchanging $x, y$ in the construction gives $\sigma_x$. 

\end{say}

\begin{say}[Proof of Theorem~\ref{Sg.main.thm}]\label{Sg.main.thm.pf}

If $a_n, b_m$ and $a_0=b_0$
are units in $R$, then
the composition of the  isomorphisms  $\sigma_x, \sigma_y$
constructed in (\ref{gen.not.alg.say})
gives  the automorphism
        $$
\sigma_{a,b}: S_{a,b}\stackrel{\sigma_x}{\longrightarrow} S_{\bar a,b}\stackrel{\sigma_y}{\longrightarrow} S_{\bar a,\bar b}\stackrel{\sigma_x}{\longrightarrow} S_{a,\bar b}\stackrel{\sigma_y}{\longrightarrow} S_{a,b}.
$$
We need to show that it has infinite order. For this it is enough to show that
the birational maps
$$
\sigma'_x: (x,y)\mapsto \bigl(x, a(x)y^{-1}\bigr)\qtq{and}
\sigma'_y: (x,y)\mapsto \bigl(b(y)x^{-1}, y\bigr)
$$
generate an  infinite  subgroup of $\bir(\a^2)$.
Now observe that 
if $\deg a, \deg b\geq 3$ and $\max\{|x|,|y|\}$ is large,  then
either $\sigma'_x$ or $\sigma'_y$ increases  $\max\{|x|,|y|\}$.
\qed

  \end{say}

        \section{Isomorphisms}

We are now ready to describe all isomorphisms in the cubic case.

        \begin{thm} \label{main.isom.thm} Let $R$ be an integral domain,  and consider the set of all cubic surfaces
          $\{ S_{a,b}: a, b\in P^\times_3(R)\}$   with $P^\times_3$ as in  (\ref{a(t).not}) and (\ref{Sab.def.R}).
          
          Then the groupoid of isomorphisms between the surfaces
          $ S_{a,b}$ is generated by
          \begin{enumerate}
          \item linear isomorphisms (\ref{norm.f.lem}.2),
          \item $\sigma_y: S_{a,b}\cong S_{\bar a, b}$
            (\ref{gen.not.alg.say}) for all $a, b\in P^\times_3(R)$, and
          \item $\sigma_x: S_{a,b}\cong S_{ a, \bar b}$
             (\ref{gen.not.alg.say}) for all $a, b\in P^\times_3(R)$.
              \end{enumerate}
        \end{thm}

Next we give some consequences  over $\z$.
        
        \begin{notation}\label{gen.aut.cor.exmp}
        Assume now that we work over $\z$ and let  $a, b\in P^\times_3(\z)$. We can  assume that $a_0=b_0=1$.
        By changing the sign of $x,y,z$ we can also assume that $a_3=b_3=1$.
        Thus we have
        $$
        \begin{array}{l}
          a(t):=t^3+a_2t^2+a_1t+1, \quad b(t):=t^3+b_2t^2+b_1t+1,\qtq{and}\\
          \bar a(t):=t^3+a_1t^2+a_2t+1, \quad \bar b(t):=t^3+b_1t^2+b_2t+1.
          \end{array}
        $$
  There are 8 surfaces in play
$$
S_{a,b},S_{\bar a,b},S_{a,\bar b},S_{\bar a,\bar b},
S_{b,a},S_{\bar b,a},S_{b,\bar a},S_{\bar b,\bar a}.
\eqno{(\ref{gen.aut.cor.exmp}.1)}
$$
We always have the obvious linear isomorphism 
$\tau: S_{a',b'}\cong  S_{b',a'}$,  so
the first 4 surfaces are mapped by $\tau$ to the last 4.

With these choices, $ \sigma_y: S_{a,b}\cong S_{\bar a,b}$ is given by 
         $$
\bigl(yz-(x^2+a_2x+a_1), y, yz^2-z(x^2+a_2x+a_1)-(x+a_2)(y^2+b_2y+b_1)\bigr),   
\eqno{(\ref{gen.aut.cor.exmp}.2)}
$$
and  $ \sigma_x: S_{a,b}\cong S_{a,\bar b}$ is given by
$$
\bigl(x, xz-(y^2+b_2y+b_1),  xz^2-z(y^2+b_2y+b_1)-(y+b_2)(x^2+a_2x+a_1)\bigr).      
\eqno{(\ref{gen.aut.cor.exmp}.3)}
$$
        \end{notation}

         Using the above notation and assumptions, we have 3 corollaries.

        \begin{cor}  \label{gen.aut.cor.1}  If  $S_{a,b}\cong  S_{a',b'}$ for some  $a', b'\in P^\times_3(\z)$, then  $ S_{a',b'}$ is one  of the 8  surfaces in (\ref{gen.aut.cor.exmp}.1). \qed
          \end{cor}

        \begin{cor}  \label{gen.aut.cor.3} If   $a(t)=b(t)=t^3+1$, then
          the 8  surfaces in (\ref{gen.aut.cor.exmp}.1) agree with $S_0$, and 
          $\aut_{\z} (S_0)=\langle  \sigma_x, \sigma_y, \tau \rangle$.
         \qed
          \end{cor}

\begin{cor}  \label{gen.aut.cor.2} Assume in addition that $a_2\neq a_1, b_2\neq b_1$ and $(a_2, a_1)\neq (b_2, b_1)\neq  (a_1, a_2)$. Then
          \begin{enumerate}
                           \item all linear isomorphisms between  the  surfaces (\ref{gen.aut.cor.exmp}.1)  are given by $\tau$, and 
                           \item  
 $\aut_\z( S_{a,b})=\aut_\c( S_{a,b})$ is infinite, cyclic, and  generated by the composition
        $$
        \sigma_{a,b}: S_{a,b}\stackrel{\sigma_x}{\longrightarrow} S_{\bar a,b}\stackrel{\sigma_y}{\longrightarrow} S_{\bar a,\bar b}\stackrel{\sigma_x}{\longrightarrow} S_{a,\bar b}\stackrel{\sigma_y}{\longrightarrow} S_{a,b}. \hfill \qed
        $$
        \end{enumerate}
\end{cor}

{\it Remark \ref{gen.aut.cor.2}.3.} One should think of $\aut_\z( S_{a,b})\cong \z$
as an index 8 subgroup of $\aut_\z( S_0)\cong D_{\infty}$.

\medskip
{\it Remark \ref{gen.aut.cor.2}.4.}
We computed the coordinate functions of $\sigma_{a,b}$ using Macaulay2.
They are polynomials of degrees $13$,  $34$, and  $55$, having 
178, 3485  and  15314  monomials.
Simplifying the formulas using a Gr\"obner basis,
the degrees stay the same, but the number of  monomials drops to
110, 998, and  2881.
We also computed some very special cases, for example
$(a_2, a_1, b_2, b_1)=(0,1,0,2)$ or  $(a_2, a_1, b_2, b_1)=(0,1,0,1)$.
The number of  monomials stays the same.

\begin{say}[Plan of the proof of Theorem~\ref{main.isom.thm}]  \label{main.isom.thm.plan}
  We try to follow the proof in \cite{k-lk3}. 

  Given any isomorphism $\psi:S_{a,b}\cong S_{a',b'}$,
  we want to  show that there is a $\phi\in \langle \sigma_x, \sigma_y\rangle$
  such that the composite
  $$
  \phi\circ \psi: S_{a,b}\cong S_{a',b'}\cong S_{a'',b''} \qtq{is linear.}
  \eqno{(\ref{main.isom.thm.plan}.1)}
  $$
  Now note that being a linear isomorphism does not depend on the base ring, and
   $\sigma_x, \sigma_y$
  are defined over $R$. Thus once we prove 
  (\ref{main.isom.thm.plan}.1) over the quotient field of $R$, it holds over $R$ as well. We may even replace the quotient field by its algebraic closure.

Thus assume from now on  that $R=K$ is an algebraically closed field, and
let $A\subset \bar S_{a,b}$ be a smooth hyperplane section with only 1 place at infinity as in (\ref{flex.tg.1.p.exmp}.1).

  Given any isomorphism $\psi:S_{a,b}\cong S_{a',b'}$, we consider the curve
  $\psi_*(A)\subset \bar S_{a',b'}$. It has only 1 place at infinity.

  Next we use isomorphisms $\sigma_x, \sigma_y$ to simplify the singularity of any curve $C$ that  has only 1 place at infinity.
 A similar question was considered in  \cite[7]{k-lk3}.  Although \cite[7]{k-lk3} applies only to smooth Del~Pezzo surfaces of degree $>3$,  extrapolating from \cite[7.1]{k-lk3}
 suggests that 
  $$
 \mult_s\bigl(\phi_*(C)\bigr)\mbox{ might be  at most } \tfrac13 \deg \bigl(\phi_*(C)\bigr).
 \eqno{(\ref{main.isom.thm.plan}.2)}
  $$
    The presence of the $A_2$ singularity changes the computation.
  The first part goes even better, see (\ref{mult.red.prop}).
  Applying it to $C:=\psi_*(A)$
  we can achieve that
  $$
  \mult_s(\phi\circ \psi)_*(A){\leq} \tfrac14 \deg\bigl((\phi\circ \psi)_*(A)\bigr).
   \eqno{(\ref{main.isom.thm.plan}.3)}
  $$
  On the other hand, the curve $(\phi\circ \psi)_*(A)$ still passes through the
  $A_2$ singularity, hence the pair
  $\bigl(\bar  S_{a'',b''}, \epsilon\cdot (\phi\circ \psi)_*(A)\bigr)$ is not canonical (\ref{discrep.say})  for
  any $\epsilon>0$.  This is bad news since the Noether-Fano method (Section~\ref{n-f.meth.sec})  would need
  $\bigl(\bar  S_{a'',b''}, \epsilon\cdot (\phi\circ \psi)_*(A)\bigr)$ to be canonical for
$\epsilon_0:=3/\deg \bigl((\phi\circ \psi)_*(A)\bigr)$.

  Instead, we are only able to obtain the  discrepancy (\ref{discrep.say}) inequality
  $$
  \discrep\bigl(\bar  S_{a'',b''}, \epsilon_0\cdot (\phi\circ \psi)_*(A)\bigr)\geq -\tfrac12.
    \eqno{(\ref{main.isom.thm.plan}.4)}
    $$
    (Canonical would be $\discrep\geq 0$.)
    The traditional Noether-Fano method needs strong control of the singularities of $(\phi\circ \psi)_*(A)$, but makes no assumption about the indeterminacy locus of $\phi\circ \psi$. In the present setting we know that the indeterminacy locus is contained in the boundary $\bar B$, and we are able to use this to relax the condition on the singularities. A result of this type is proved in (\ref{n-f.meth.thm}).

     We put the pieces together in (\ref{main.isom.thm.pf}). The end result is that $(\phi\circ \psi)_*(A)$
     lies either in $|-K|$ or in $|-2K|$.
  In the first case, $\phi\circ \psi$ is linear, as needed. The second case is  excluded using (\ref{l.les.q.lem}.3). 
 \end{say}

  \begin{defn} \label{flex.tg.1.p.exmp}
 Let $C$ be an irreducible curve, with
  projective closure $\bar C$. We say that $C$ has
  {\it one place at infinity} if the normalization of  $\bar C$
  has only one  geometric point  over  $\bar C\setminus C$. Equivalently,
  $\bar C$ has a  single geometric  point at infinity, and it is a cusp.

  {\it Example \ref{flex.tg.1.p.exmp}.1.} Let $\bar S$ be a cubic surface whose curve at infinity $\bar B$ is nodal.
  Let $p\in \bar B$ be a flex with flex tangent line  $L_p$.
  For a general plane $H\supset L_p$, the intersection
  $A:= S\cap H$ is a smooth, elliptic curve with one place at infinity.

  The  nodal cubic $\bar B$ need not have a flex over the base field $k$, 
but it
has 3 flexes over the algebraic closure. So there is always a flex over a suitable cubic extension of $k$.
\end{defn}

\section{Transforming curves}\label{tra.cur.sec}

\begin{notation} Let $k$ be a field and $S=S_{a,b}$  an affine  surface as in (\ref{Sab.def.R}). Let  $\bar S$ be its projective closure,  $\bar s\in \bar S$ the unique singular point at infinity,  and  $\bar B\subset \bar S$ the curve at infinity. Let $\bar T\to \bar S$ denote the minimal resolution, with exceptional curve $E=E_x\cup E_y$.

  If $C$ has
  one place at infinity,
  we use $\mult_E(C)$ to denote the multiplicity of  the birational transform  of $C$ on $\bar T$ at its  intersection point with $E$.
      Note that $\mult_E(C)=0$ if $\bar s\not\in \bar C$. Also, 
    $\mult_{\bar s}\bar C\geq\mult_E C$, but they need not be equal.

\end{notation}

The following is a key step toward proving (\ref{main.isom.thm}).

\begin{prop} \label{mult.red.prop} Let $C\subset S$ be a curve with one place at infinity. Then there is a $\phi\in \langle\sigma_x, \sigma_y\rangle$
  giving
  $\phi: (S, C)\map (S^{\rm m}, C^{\rm m})$ such that either $\deg C^{\rm m}<3$,
  or 
  $\mult_E C^{\rm m}\leq \tfrac14\deg C^{\rm m}$.
\end{prop}

Proof. The transformations $\sigma_x, \sigma_y$ both involve
2 blow-ups followed by 2 contractions. This is exactly the $r=2$ special case of
\cite[29]{k-lk3}. Using its notation,  $p=\mult_E(C)$ and
$q$ is the intersection multiplicity of the birational transform  of $C$ with  the birational transform  of $\bar B$. So $\deg C=p+q$. 
We check in  (\ref{l.les.q.lem}) that $p<q$. 
Thus
$B\stackrel{p,q}{{\mbox{ --- }}} B$ is transformed to
 $$
    \begin{array}{ll}
      B \stackrel{p, 2p-q}{{\mbox{ --- }}} B
       &
      \mbox{if $q\leq 2p$, and}\\[1ex]
       B \stackrel{q-2p,p+2(q-2p)}{{\mbox{ --- }}} B
       &
      \mbox{if $q\geq 2p$.}
    \end{array}
    \eqno{(\ref{mult.red.prop}.1)}
    $$
Equivalently, 
$$
\deg\bigl(\sigma_x(C)\bigr)=
    \left\{
      \begin{array}{l}
        4\mult_E(C) -\deg(C)\qtq{if \ $\mult_E (C)\geq \tfrac13\deg (C)$, and}\\[1ex]
        3\deg(C)-8\mult_E(C)\qtq{otherwise.}
      \end{array}
      \right.
      $$
      In the first case, $4\mult_E(C) -\deg(C)< \deg(C)$ iff
      $\mult_E(C) <\tfrac12\deg(C)$. This holds by
 (\ref{l.les.q.lem}.2), unless  $\deg C<3$.

      In the second case, $3\deg(C)-8\mult_E(C)< \deg(C)$ iff
      $\mult_E(C) > \tfrac14\deg(C)$.

      We can thus lower the degree using the transformations
      $\sigma_x, \sigma_y$, except when we have reached
      $(S^{\rm m}, C^{\rm m})$ such that either $\mult_E C^{\rm m}\leq \tfrac14\deg C^{\rm m}$ or $\deg C^{\rm m}<3$.\qed

      \begin{lem}\label{mult.bd.1.lem} Let $T$ be a projective,   weak Del~Pezzo surface  of degree 3, and $C\subset T$ an irreducible curve that is not
        contained in an effective anticanonical divisor. Then $\mult_cC\leq \tfrac12 \deg C$ for every smooth point $c\in T$.
      \end{lem}

      Proof.  Since $\dim |-K_T|=3$, there is an anticanonical divisor $F_c$ that is singular at $c$. By assumption $C\not\subset F_c$. Thus
      $$
      2\mult_cC\leq \mult_cF_c\cdot \mult_cC\leq (F_c\cdot C)=\deg C. \qed
      $$

\begin{lem}[Curves with 1 place at infinity]\label{l.les.q.lem}
  Let $C\subset \bar S$ be a curve of degree $\geq 3$ with 1 place at infinity at the node.
  Let $\pi:\bar T\to \bar S$ be the minimal resolution with exceptional curves $E_1, E_2$. Let $C_T$ be the birational transform of $C$ on $\bar T$. Then
  \begin{enumerate}
  \item $C_T$ meets $E_1\cap E_2$ at a point of $\bar B_T$, 
  \item  $\bigl(C_T\cdot(E_1+ E_2)\bigr)< \bigl(C_T\cdot \bar B_T\bigr)$, and
    \item there is no such curve in $|-2K_{\bar S}|$.
  \end{enumerate}
  \end{lem}

Proof.  If (1) does not hold, then $C_T$ is disjoint from $\bar B_T$, which is the pull-back of a line in the blow-up model (\ref{blw.up.say}). Thus $C_T$ is one of the 6 exceptional curves of the blow-up model, so $C$ is a line through the singular point.

Assume next that the common point is $\bar s:=E_2\cap \bar B_T$.
Set $p:=\bigl(C_T\cdot(E_1+ E_2)\bigr)$ and $q:=\bigl(C_T\cdot \bar B_T\bigr)$. 
 Blow-up $\bar s$. Denote the exceptional curve by $E_3$, and let
 $C_3, B_3$ be   the birational transforms $C_T, \bar B_T$.

 If  $p\geq q $ 
then $C_3$ is disjoint from $B_3$ and has intersection number $p$ with $E_3$.
Note that $(B_3^2)=0$, so it moves in a basepoint-free pencil, with section $E_3$.
Then $C_3$ is a subcurve of a member of the pencil, hence $(C_3\cdot E_3)\leq 1$.
So $p=1$ and $\deg C\leq 2$.
These cases are listed in (\ref{blw.up.say}).

If $C\in |-2K_{\bar S}|$ , then $C$ is Cartier, so  $\pi^*C=C_T+\frac{2p}{3}E_2+\frac{p}{3}E_1$ and
$\deg C=\bigl(C_T\cdot(E_1+ E_2)\bigr)+ \bigl(C_T\cdot \bar B_T\bigr)$.
If $\deg C\leq 6$ then $p<3$ by (2), hence $\frac{p}{3}$ is not an integer. 
 \qed

 \begin{cor} \label{unsplit.node.cor}
   Let $k$ be a field and $S:=\bigl(zg_2(x,y)=g_3(x,y)\bigr)$  an affine  cubic surface as in  (\ref{our.cubics.defn}). 
   If $g_2(x, y)$ is irreducible over $k$, then 
   every $k$-automorphism of $S$ is linear.
 \end{cor}

 Proof. By (\ref{flex.tg.1.p.exmp}), after possibly a degree 3 extension $k'/k$, there is a curve
 $C_1\subset S$ with  1 place at infinity. Then $g_2(x, y)$ is still irreducible over $k'$.
 If there is a nonlinear $k'$-automorphism, then the image of $C_1$ is a
 curve  $C\subset \bar S$  of degree $\geq 3$ with 1 place at infinity at the node. By (\ref{l.les.q.lem}.1), then  $C_T$ meets $E_1\cap E_2$ at a $k'$ point of $\bar B_T$. If the node is unsplit, this intersection is a degree 2 point over $k'$. \qed

 \section{Noether-Fano method}\label{n-f.meth.sec}

 We start with some general comments on  closed graphs of birational maps.
     
\begin{say}\label{must.node.lem}
  Let  $T_i$ be proper surfaces  with rational singularities and
  $D_i\subset T_i$   geometrically irreducible  curves.
Let   $\phi^\circ: T_1\setminus D_1 \cong T_2\setminus D_2$ be an isomorphism
   that 
 does not extend to an isomorphism of 
 $T_1$ and $T_2$.

 Let $\bar T$ be the normalization of the closed graph of $\phi^\circ$ with projections $p_i: \bar T\to T_i$.  Let $\bar D_i$ be the birational transform of $D_i$ on $\bar T$.

 Let $E_i\subset \bar T$ be a $p_i$-exceptional curve. Then $E_i$ is not
 $p_{3-i}$-exceptional, hence  $p_{3-i}$ must map $E_i$ birationally onto $D_{3-i}$.
 Thus  $E_i=\bar D_{3-i}$ is the unique exceptional curve of $p_{i}$.
 
 Since the $T_i$ have rational singularities, the  $E_i$ are smooth rational curves, contracted to a point
 $t_i\in T_i$.
 Since $\bar D_i\to D_i$ is an isomorphism away from $t_i$, we conclude that
 $t_i$ is the only possible singular point of $D_i$.

 Assume now that the $D_i$ are nodal and the pairs  $(T_i, D_i)$ are log canonical. As in
 \cite[17]{k-lk3} we get that $\bar D_1, \bar D_2$ meet at 2 points on $\bar T$,
 both cyclic quotients. On the minimal resolution of $\bar T$, the exceptional curves and the  birational transforms of the $D_i$ form a cycle of rational curves.
 
  \end{say}

  These are strong restrictions, but still leave many possibilities.
 In order to exploit (\ref{mult.red.prop}),
we need a few definitions and results on discrepancies.
  We state these for surfaces, which is the only case that we use.
     For general introductions, see
     \cite[Sec.2.3]{km-book} or \cite[Sec.2.1]{kk-singbook}.

  \begin{say}[Discrepancies for surface pairs]\label{discrep.say}
    Let $S$ be a normal surface and $\pi:T\to S$ a proper, birational morphism.
We are interested in the local behavior of $\pi$ over a point $s\in S$, 
  so   assume that $\pi$ is an isomorphism over $S\setminus\{s\}$, and $T$ is normal.
    We can write $K_T\simq \pi^*K_S+E$ where $E$ is  $\pi$-exceptional (with rational coefficients).

    Next let $\Delta=\sum a_iD_i$ be a finite, $\q$-linear combination of  distinct curves on $S$, and $\Delta_T$ its birational transform on $T$. Then  $\Delta_T\simq \pi^*\Delta-F$ where  $F$ is  $\pi$-exceptional.  We can formally write
    $$
    K_T+\Delta_T\simq \pi^*(K_S+\Delta)+E-F.
    \eqno{(\ref{discrep.say}.1)}
    $$
    The coefficient of an exceptional curve $E_j$  in $E-F$ is called the {\it discrepancy} of $E_j$,  denoted by  
    $a(E_j, S, \Delta)$.  We set $a(D_i, S, \Delta):=-a_i$.

    The infimum over all exceptional curves and all
    proper, birational morphisms is called the
    {\it discrepancy} of 
    $(S, \Delta)$, denoted by   $\discrep(S, \Delta)$.

    If $\Delta$ is effective, 
    the pair $(S, \Delta)$ is called  {\it canonical} if  $\discrep(S, \Delta)\geq 0$, and {\it log canonical} if  $\discrep(S, \Delta)\geq -1$.

    We need 2 results.
    The first is a special case of \cite[2.30]{km-book}. The second
    follows from  a direct combination of  \cite[5.50]{km-book}  and \cite[2.5]{kk-singbook}. For surfaces it can be proved directly by computing one blow-up, and using induction.

    \medskip
    {\it Claim \ref{discrep.say}.2.}
    $
    \discrep(S, \Delta)=\min\bigl\{ \discrep(T, \Delta_T+F-E), a(E_j, S, \Delta)\bigr\},
    $
    where $E_j$ runs through all  $\pi$-exceptional curves.\qed

      \medskip
    {\it Claim \ref{discrep.say}.3.}  Let $S$ be a smooth surface,
    $D_0\subset S$ a smooth curve, and $\Delta$ an effective $\q$-divisor such that
    $(D_0\cdot \Delta)\leq 1$. Then, in a neighborhood of $D_0$, 
    $    \discrep(S, d_0D_0+\Delta)\geq -d_0$. \qed

  \end{say}

  We can now translate (\ref{mult.red.prop}) into a discrepancy statement.
  First a general version.

\begin{lem} \label{1/3.lc.lem}
  Let $(t\in T)$ be an $A_2$ singularity, and $\pi:T'\to T$ its minimal resolution with exceptional curves  $E_1, E_2\subset T'$. Let $C\subset T$ be a curve and $C'$ its birational transform on $T'$. Assume that
  $(E_1\cdot C')=0$,  $(E_2\cdot C')=r$, and $C'$ has a cusp of multiplicity $r$ at $E_2\cap C'$. Then, for $0\leq c\leq \tfrac1{r}$, the  discrepancy of 
  $(T, cC)$ is $-2cr/3$ (in a neighborhood of $t$).
\end{lem}

Proof.  First note that
$c\pi^*C\simq  cC'+(2cr/3) E_2+(cr/3) E_1$.
Using (\ref{discrep.say}.2) it remains to show that  discrepancy of  $\bigl(T', cC'+(2cr/3) E_2+(cr/3) E_1\bigr)$  is at least $-2cr/3$.
This follows from (\ref{discrep.say}.3).\qed

\medskip

\begin{cor} \label{mult.red.prop.cor} Let $C^m\subset S^m$ be
as in (\ref{mult.red.prop}), and assume that $C^m$ is smooth away from  $\sing \bar S^m$. 
  Then the discrepancy  of
  $\bigl(\bar S^m, \tfrac3{\deg C^m} C^m\bigr)$ is at least $-\tfrac12$.
\end{cor}

Proof. Set  $r=\mult_E C$ and apply (\ref{1/3.lc.lem}) with $c=\tfrac3{\deg C^m}$. \qed
    
\medskip
  
In order to study automorphisms of open varieties, we need a variant of the 
Noether-Fano method. 
See  \cite[Sec.5.1]{ksc} for a description of the classical version, and \cite[Sec.8]{k-lk3} for related results.
We state the general case, but then use it only for surfaces.

\begin{thm}\label{n-f.meth.thm}
  Let $X_i$ be proper, normal varieties, $D_i\in  |-K_{X_i}|$  irreducible divisors and  $\phi^\circ:(X_1\setminus D_1)\cong (X_2\setminus D_2)$  an isomorphism.
  Let $D_1\neq A_1\in  |-K_{X_1}|$ be a divisor
  and  set $A_2:=\phi^\circ_*(A_1)$. Assume that
  $A_2\in  |-mK_{X_2}|$ and the discrepancy of $\bigl(X_2, \tfrac1{m} A_2\bigr)$ is at least 
  $-1+\epsilon$.
Then $m\epsilon\leq 1$.
\end{thm}

{\it Comment \ref{n-f.meth.thm}.1.} The traditional Noether-Fano method does not involve the divisors $D_i$, but
assumes that $\bigl(X_2, \tfrac1{m} A_2\bigr)$ is canonical, that is,
$\epsilon=1$. It gives that  $m=1$, hence we get an isomorphism $X_1\cong X_2$, provided the $-K_{X_i}$ are ample.

In our case the $D_i$ put strong restrictions on the indeterminacy loci of $\phi$ and $\phi^{-1}$, but we have weaker assumptions on the singularities of $\bigl(X_2, \tfrac1{m} A_2\bigr)$.
\medskip

Proof.   Let $Y$ be the normalization of the closure of the graph of
$\phi^\circ$ with projections $p_i:Y\to X_i$.
Let  $A_Y$ be the birational transform of the $A_i$ on $Y$.  Write
$K_Y\simq p_i^*K_{X_i}+E_i$ and $A_Y=p_i^*A_i-F_i$.
For any $c\in \q$  we have
 $$
K_Y+cA_Y\simq    p_i^*\bigl(K_{X_i}+cA_i\bigr)+E_i-cF_i.
\eqno{(\ref{n-f.meth.thm}.2)}
 $$
Now choose  $c=\frac1{m}$. Then $K_{X_2}+\tfrac1{m}A_2\simq 0$ and
$K_{X_1}+\tfrac1{m}A_1\simq \tfrac{m-1}{m}K_{X_1}$, hence
  $$
  E_2-\tfrac1{m}F_2 \simq 
 p_1^*\bigl(\tfrac{m-1}{m}K_{X_1}\bigr)+E_1-\tfrac1{m}F_1.
 \eqno{(\ref{n-f.meth.thm}.3)}
 $$
 If  $p_2\circ p_i^{-1}$  is an isomorphism at the generic points of the $D_i$, then linear equivalence is preserved, hence $m=1$. 
 Otherwise $D_1$ is the image of
 an irreducible component  $D_1^Y$ of $E_2-\tfrac1{m}F_2$, whose
 coefficient is $d_1\geq -1+\epsilon$
 by assumption.  Thus
  $$
 d_1D_1=(p_1)_*\bigl( E_2-\tfrac1{m}F_2\bigr)\simq \tfrac{m-1}{m}K_{X_1}
 \simq -\tfrac{m-1}{m}D_1.
 $$
 Hence $-\tfrac{m-1}{m}=d_1\geq -1+\epsilon$, and so  $m\epsilon\leq 1$.
 \qed

 \begin{cor} \label{main.isom.k.cor}
   Let $A\subset \bar S_{a,b}$ be a smooth hyperplane section with only 1 place at infinity, and $\psi:S_{a,b}\cong S_{a',b'}$
   an isomorphism. Then  there is a
    $\phi\in \langle \sigma_x, \sigma_y\rangle$
  such that under the composite
  $$
  \phi\circ \psi: S_{a,b}\cong S_{a',b'}\cong S_{a'',b''}, 
  \eqno{(\ref{main.isom.k.cor}.1)}
  $$
  the image  $(\phi\circ \psi)_*(A)$ lies in $|-K|$ or $|-2K|$ on $S_{a'',b''}$.
 \end{cor}

 Proof. Apply (\ref{mult.red.prop}) and (\ref{mult.red.prop.cor}) to $\psi_*(A)$.
 We get a  $\phi\in \langle \sigma_x, \sigma_y\rangle$
 such that
 $$
 \bigl( S_{a'',b''}, \tfrac1{m} (\phi\circ \psi)_*(A)\bigr)
 $$
 has discrepancy $\geq -\tfrac12$, where
 $m=\frac13\deg\bigl((\phi\circ \psi)_*(A)\bigr)$.

 As we noted in \cite[16]{k-lk3}, the class group of
 $S$ is $\cl(\bar S)/\z[-K_{\bar S}]$. Thus the class of $A$ is trivial in
 $\cl(S_{a,b})$, so
  the class of $(\phi\circ \psi)_*(A)$ is trivial in
 $\cl(S_{a'',b''})$. Therefore 
  $(\phi\circ \psi)_*(A)\in |-mK|$ for some $m$ of $S_{a'',b''}$.

  We can now apply 
  (\ref{n-f.meth.thm}) to get  that $m\leq 2$. \qed
 
 \begin{say}[End of the proof of Theorem~\ref{main.isom.thm}]
   \label{main.isom.thm.pf}
Using (\ref{main.isom.k.cor}), we get
    $\phi\in \langle \sigma_x, \sigma_y\rangle$
  such that 
  $A^m:=(\phi\circ \psi)_*(A)$  lies in $|-K|$ or $|-2K|$  on $S_{a'',b''}$.
  
The possibility $A^m\in |-2K|$ is excluded by
(\ref{l.les.q.lem}.3). Thus
 $A^m\in |-K|$ and so  $\phi\circ \psi$ is a linear isomorphism. \qed
 \end{say}

 \section{Integral points}\label{int.pts.sec}

We do not have a good   understanding of all integral points on cubic surfaces, but  Zariski density is known in many cases.

 \begin{say}[Cubics containing a line]\label{cubic.line.say}
   Let $K$ be a number field with ring of integers $R$.
   Let $S$ be an affine, cubic surface that contains a $K$-line $L$.
    Projecting $S$ from $L$ gives a map $\pi:S\setminus L\map \p^1$ whose general fibers are $\gm$-torsors. This is not very helpful if $\bar S$ has a singular point on $\bar L$. However, if $\bar S$ is smooth along $\bar L$, then $\pi$
   was used to prove the Zariski density of $R$-points in increasing generality in the papers
   \cite{MR51852, MR78397, MR1352640, MR1684588, MR1875174}; see also \cite{MR3908761} for related results.
   The surface $S_0$ is treated already in  \cite{MR51852}.

   As we noted in (\ref{lines.Sab.say}), 
         there is such a  $K$-line on $S_{a,b}$ iff both  $a(x)$ and $b(y)$ have a root in $K$. Next we see how to use this approach to get integral points.
   \end{say}

 \begin{say}\label{get.log2.gr.say}
   Let $S_{a,b}$ be as in (\ref{norm.f.lem}.1).
   As in (\ref{blw.up.say}), we think of it as  the blow-up of $\a^2$ at 6 base points, corresponding to the roots of  $a(x)$ and  $b(y)$.
   
   Let $\alpha$ be a root of $a(x)$ and
    $\beta$ a root of $b(y)$.  Write 
   $
   a(x)=(\alpha^{-1}x-1)p(x)$ and  $  b(y)=(\beta^{-1}y-1)q(y)$.
   Then
      $$
   Q_{\lambda,\mu}:=\lambda \bigl(p(x)+q(y)-c\bigr)+\mu xy
   $$
   is a pencil of quadratic polynomials, vanishing at  4 of the 6 base points.
   Thus on $S_{a,b}$ the birational transforms $C_{\lambda,\mu}$ of the curves $(Q_{\lambda,\mu}=0)$ form a pencil of conics residual to the line $L_{\alpha, \beta}$, which is the birational transform of $\bigl(\alpha^{-1}x+\beta^{-1}y=1\bigr)$.

   Assume now that $\alpha, \beta$ are rational and 
   let $(x_0, y_0, z_0)$ be an  integral point on $S_{a,b}$.
   There is a unique (up to scalar) $Q_{\lambda_0,\mu_0}$ that vanishes at
   $(x_0, y_0)$.

   If we are lucky, then $Q_{\lambda_0,\mu_0}(x,y)=0$ has infinitely many integer solutions. This is not automatic, and we used \cite{alpern24} to check cases.

   If there are infinitely many integer solutions, then they are described by starting values and a linear recursion. This implies that
   \begin{enumerate}
   \item the number of solutions of $Q_{\lambda_0,\mu_0}(x,y)=0$ grows  like
     $\epsilon\log N$, and
   \item  the solutions  are periodic modulo  any integer.
     \end{enumerate}

   We need to understand when the corresponding $z$ values are integers.
 The $C_{\lambda,\mu}$ have 2 points at infinity, hence 
   the restriction of $z$ to  $C_{\lambda,\mu}$ is a linear function in $x,y$, with rational coefficients. Thus (\ref{get.log2.gr.say}.2) implies that $z$ is 
   periodic modulo 1.
   So, if there is one value  $(x_0, y_0)$ for which $z_0$ is an integer,
   then there are infinitely many.
   Also, in these cases the  $Q_{\lambda_0,\mu_0}=0$ are hyperbolas, hence  one of the branches lies in the fundamental domain
   for $x\gg 1$.

   The automorphism group then moves these solutions to get new ones, but these seem to grow double exponentially. 
   We have thus proved the following.

   \medskip
       {\it \ref{get.log2.gr.say}.3 Corollary.} If there is such a
       $Q_{\lambda_0,\mu_0}(x,y)$ and $\aut(S_{a,b})$ is infinite, then 
 the number of integral points on $S_{a,b}$ grows at least like $\epsilon\log N\log\log N$. \qed
   \medskip

   Now consider $S_0=(xyz=x^3+y^3+1)$. The rational roots are $\alpha=\beta=-1$, and
   $$
   Q_{\lambda,\mu}=\lambda (x^2+y^2-x-y+1)+\mu xy.
   $$
   For the first 3 solutions $(1,1,3) (1,2,5),  (2,3,6)$ the resulting degree 2 equation    $Q_{\lambda_0,\mu_0}=0$ has only finitely  many integer solutions.
   For the next case   $(2,9,41)$ the degree 2 equation is
   $$
   6x^2-25xy+6y^2-6x-6y+6=0
   \eqno{(\ref{get.log2.gr.say}.4)}
   $$
   which has infinitely many integer solutions; see \cite{alpern24}.
   Thus (\ref{get.log2.gr.say}.1) applies and we get the following.

   \medskip
       {\it \ref{get.log2.gr.say}.5 Corollary.}  The number of integral points on       $xyz=x^3+y^3+1$ 
        satisfying
  $1\leq x,y\leq N$ grows at least like $\epsilon\log N\log \log N$ for some $\epsilon>0$.\qed
        \medskip

         {\it \ref{get.log2.gr.say}.6 Remark.} 
         The recursion for the equation (\ref{get.log2.gr.say}.4) involves 13 digit coefficients. The next 2 cases $(3,14,66), 
         (5,9,19)$ also give infinitely many solutions; for these the  recursion
          involves 3 digit coefficients.
   \end{say}

 \begin{say}[Which power of $\log$?]\label{which.power.rem}

   In order to develop a heuristic, we
   work with the divisibility conditions $ x\mid b(y)$ and $y\mid a(x)$.
   The likelyhood that $x$ divides a random integer is $\tfrac1{x}$.
   Thus if $ b(y)$ and $a(x)$ are random, then we expect that the number of solutions satisfying $1\leq x, y\leq N$ is
   $\int_1^N\int_1^N\tfrac1{xy}dxdy=\log^2 N$.

   Taking into account that we also need $x\leq b(y)$ and $y\leq a(x)$,
   refines this guess  to $\bigl(1-\frac{n+m}{2nm}\bigr)\log^2 N$,
   where $m=\deg a, n=\deg b$. 

   By contrast,
the extension of the  Batyrev--Manin conjectures to K3 surfaces
(see for example \cite{vanluijk, MR3732687}) 
suggests that the exponent of $\log N$ should be  the Picard number (over $\q$).
Over a field $k$, the Picard number of  $ S_{a,b}$ equals
$$
\#(\mbox{irreducible factors of $a(x)$})+
\#(\mbox{irreducible factors of $b(y)$})-2.
$$
For various choices of $a(x), b(y)$ 
we searched for all positive integral solutions  up to $10^5$.
The results suggest a  Batyrev--Manin--type dependence.

A detailed theoretical and numerical study of integral points on cubic surfaces is given in 
\cite{browning2024integralpointscubicsurfaces}. For the surfaces
$S_{a,b}$ \cite[Conj.1.1]{browning2024integralpointscubicsurfaces}
predicts that  the exponent of $\log N$ should be  the Picard number plus 2.
This seems larger than what we found. Note, however, that in many examples the number of points fluctuates wildely for $N<10^5$, but stabilizes for $N\gg 10^5$;
see especially \cite[Fig.2]{browning2024integralpointscubicsurfaces}.

Further numerical work would be needed to get a clearer picture.

\end{say}

 \begin{say}[Density of orbits]\label{fin.orb.say} 
   Using the notation of     (\ref{Sab.def.R}) for cubic surfaces $S_{a,b}$, by direct computation we see that $\sigma_y$ increases the multiplicity of a curve germ. Thus a 1-dimensional, closed subscheme of $S_{a,b}$ can not be invariant under $\aut(S_{a,b})$.
   Therefore, every infinite orbit of $\aut(S_{a,b})$ is Zariski dense, and there are no nontrivial $\aut(S_{a,b})$-equivariant morphisms  from $S_{a,b}$ to another variety.  So   $\sigma_{a,b}$ is loxodromic.

   As an example, assume in addition that all coefficients of $a(x)$ and $b(y)$ are nonnegative. Then  $S_{a,b}\cap (\r^+)^3$ is  $\aut(S_{a,b})$-invariant, and 
    $x\geq y>0$ implies that $a(x)/y>x$. Thus  every  $\aut(S_{a,b})$ orbit on
   $S_{a,b}\cap (\r^+)^3$ is infinite.

   If $a(t), b(t)\in \z[t]$, then  $S_{a,b}$ has  a positive, integral point
   by \cite[Thm.4]{MR3320468}, so the
    positive, integral points are Zariski dense on $S_{a,b}$.

        \end{say}

 Any $\z$-morphism  $\a^1\to S_{a,b}$ gives infinitely many integral points of $S_{a,b}$; these are called {\it affine lines.}
The existence of  affine lines depends subtly on $a(x), b(y)$.
The surface  $S_0$ seems rather special.

 \begin{exmp}[Affine lines on $S_0$]\label{aff.l.say}
   For $P_0=\p^2\setminus(xyz=x^3+y^3)$, all $\c$-embeddings  $\a^1\into P_0$
   were enumerated in \cite[8]{k-lk3}.
    Since  $S_0\to P_0=\p^2\setminus(xyz=x^3+y^3)$ is unramified,
   each of these  lifts to 3  embeddings  $\a^1\into S_0$, permuted by $\mu$ (\ref{lin.aut.say}). We discuss these and their integral structures.

     (\ref{aff.l.say}.1)  If $x=0$ then $y^3=-1$. The rational solution is  $y=-1$,  giving the line
     $\ell_1: t\mapsto (0,-1, t)$. We get integral points for $t\in \z$. This line and its
   images by  $\aut_\z(S_0)$ are disjoint from $(\r^+)^3$.
   The same holds for $\ell_2: t\mapsto (-1,0,  t)$.

     (\ref{aff.l.say}.2) $(3x+3y+z=0)$.   Then
     $xyz-x^3-y^3=-(x+y)^3$.  So $(x+y)^3=-1$, giving the parametrization
     $\ell_3: t\mapsto  (t, -1-t, 3)$.  We get integral points for $t\in \z$. As before, this line and its
     images by  $\aut_\z(S_0)$ are disjoint from $(\r^+)^3$.
     
     (\ref{aff.l.say}.3) $(21x^2-22xy+21y^2-6xz-6yz+z^2=0)$.
     As a conic in $\p^2$, an explicit birational parametrization is given by
     $$
     \phi: (u{:}v)\mapsto (u^2-uv+v^2, u^2+uv+v^2, 2u^2+10v^2).
     $$
     For $u,v\in\r$, the image  $\phi(u,v)$ lies in  $(\r^+)^3$, but
     we check in (\ref{no.zpt.conic.cor})
     that this curve does not even have $\z_2$-points.
     Thus we do not get any integral points  in  $(\r^+)^3$ from this or its
     $\aut_\z(S_0)$-images.

     However, setting $v=1$, we get
     $t\mapsto (t^2-t+1, t^2+t+1, 2t^2+10)$, giving $\sqrt{2N}$  positive  integer solutions satisfying $1\leq x, y\leq N$, for the equation
     $xyz=x^3+y^3+8$. 
     \end{exmp}

\begin{lem}\label{no.zpt.conic.lem}  Let $R$ be a UFD such that either $2$ is invertible in $R$
  or  $\deg \bigl(R/pR : \f_2)$ is odd for every $p\mid 2$.   Then all 
  $R$-solutions of $(21x^2-22xy+21y^2-6xz-6yz+z^2=0)$ are
  $$
     \phi(u,v):=   (u^2-uv+v^2, u^2+uv+v^2, 10u^2+2v^2) \qtq{where} u,v\in R.
  $$
\end{lem}

Proof. One checks that $\phi$ is an isomorphism over $\q$.

Let $p$ be a prime in $R$ such that $p\mid \phi(u,v)$.
Then $p\mid 2uv$. So if $p\nmid 2$ then  either $p\mid u$ or $p\mid v$.
Then $p\mid u^2-uv+v^2$ gives that it divides the other too.

If $p\mid 2$ then   the only condition is  $p\mid u^2-uv+v^2$.
If $\deg \bigl(R/pR : \f_2)$ is odd then $x^2-x+1$ has no roots in $R/pR$, and then $p$ divides both $u$ and $v$. \qed

\begin{cor} \label{no.zpt.conic.cor}  $(21x^2-22xy+21y^2-6xz-6yz+z^2=0)$ has no $\z_2$ points on $S_0$.
\end{cor}

Proof.  Computing   $x^3+y^3-xyz$ on $\phi(u,v)$ we get   $8v^6$.  
This can not be 1. \qed
\medskip

\begin{exmp}[More affine lines]
  Let $S$ be a smooth cubic surface as in  (\ref{our.cubics.defn}).
  Arguing as in  \cite[2.3]{takahashi1996curves} we see that there are
  infinitely many  smooth, rational curves on $\bar S$
  with only 1 place at infinity.

  These can be constructed using  morphisms  $\rho: \bar S\to B_1\p^2$, the 1 point blow up of $\p^2$. Moreover, we get such  rational curves defined over a field $k$ iff there is such a $\rho$ defined over $k$.

  For the cubics $S_{a,b}$, such a $\rho$ exists iff all roots of $a(t), b(t)$ are in $k$. On $\bar S_0$, smooth rational curves defined over $\q$ are either lines or conics. 
\end{exmp}

For integral points the relevant question is the existence of $\z$-morphisms $\a^1\to S_0$; being an embedding is not important.
There are many other $\c$-morphisms $\a^1\to S_0$   \cite{MR2435425, mcduff2024sesqui},
but we  do not know any other $\z$-morphisms  $\a^1\to S_0$. Our methods do not seem to be able to prove that there are none.

\begin{exmp}  \label{nodal.to.sm.exmp}
  In the pencil of elliptic curves  \cite[46]{k-lk3}
  $$
  \lambda(xyz-x^3-y^3)+\mu(-xy^2-x^2z+yz^2)=0
  $$
  there are 3  rational curves with 1 place on $(xyz=x^3+y^3) $.
  One of them,  for $(\lambda:\mu)=(-3:1)$,  is defined over $\q$.
  It can be parametrized as
  $$
  t \mapsto (-t^2+t-1,t,t^3-2t^2+3t-3),
  $$
   giving a family of integer solutions, but no 
   positive ones.
   
     Note that, as a curve in $\p^2$, it has a node at $(2{:}1{:}5)$.
  The corresponding $t$ values are the nontrivial cube roots of 1.
  When we lift to $S_0$, the curve becomes smooth, and the point $(2,1,5)$ is not on it.
  This curve turns out to be the same as
  $\sigma_y\circ \ell_3$ (\ref{aff.l.say}.2).
  \end{exmp}

  \begin{say}[Another form of $S_{a,b}$]\label{div.cong.say}
It is easy to see that  for any $n,m$, the pair of
  divisibility conditions
  $$
  x\mid y^m+1\qtq{and} y\mid x^n+1
  $$
  is equivalent to the equation
  $  xyz=x^n+y^m+1.$
  More generally, set
    $$
    S'_{a,b}:=\bigl(yu-a(x)=xv-b(y)=0\bigr)\subset \a^4_{xyuv}.
    \eqno{(\ref{div.cong.say}.1)}
    $$
There is a natural morphism $S_{a,b}\to S'_{a,b}$ given by
    $ u:=xz-b^*(y)$ and $v:=yz-a^*(x)$.
    Then $cz=uv-a^*(x)b^*(y)$, hence $S_{a,b}\cong S'_{a,b}$ over $R$ if $c$
    is a unit.

    If $a(x), b(y)$ are cubics, then the closure of $ S'_{a,b}$ in $\p^4$ is a non-normal surface whose canonical class is ample, so, for geometric questions  it is harder to work with. On the other hand,  $ S'_{a,b}$
    exists even if $a(0)\neq b(0)$, giving different integral structures on  the
    surfaces  $ S_{a,b}$.

    Furthermore, for any polynomials  $a_i(x), b_i(y)$ there is a natural morphism
    $$
    S'_{a_1, b_1}\to  S'_{a_1a_2, b_1b_2}.
    \eqno{(\ref{div.cong.say}.2)}
    $$
    Thus, if $a_1(0)=b_1(0)=a_2(0)=b_2(0)=1$, then we get
     a natural morphism
    $$
     S_{a_1, b_1}\to  S_{a_1a_2, b_1b_2}.
     \eqno{(\ref{div.cong.say}.3)}
    $$
  \end{say}

  \begin{exmp}\label{div.cong.say.exmp}
  Applying (\ref{div.cong.say}.3)  to $S_0=(xyz=x^3+y^3+1)$ we get
  4 surfaces with $\z$-morphisms to it. These are
    $(xyz=x^2+y^2-x-y+1)$,
   $(xyz=x^2-x+y+1)$,
   $(xyz=y^2-y+x+1)$, and
   $(xyz=x+y+1)$.
    It is easy to see that  the integral points are not Zariski dense on the last 3 surfaces. For example, for $(xyz=x^2-x+y+1)$
  write  $xu=y+1$. Then $xu-1\mid x^2-x+1$ and 
  $xu-1\mid ux^2-ux+u=x(xu-1)-(xu-1)+x+u-1$.
   Thus $xu-1\mid x+u-1$.

   The  integral points on $(xyz=x^2+y^2-x-y+1)$ are also not Zariski dense;
   see \cite[pp.299-300]{MR249355}.
     \end{exmp}



  \def\cprime{$'$} \def\cprime{$'$} \def\cprime{$'$} \def\cprime{$'$}
  \def\cprime{$'$} \def\dbar{\leavevmode\hbox to 0pt{\hskip.2ex
  \accent"16\hss}d} \def\cprime{$'$} \def\cprime{$'$}
  \def\polhk#1{\setbox0=\hbox{#1}{\ooalign{\hidewidth
  \lower1.5ex\hbox{`}\hidewidth\crcr\unhbox0}}} \def\cprime{$'$}
  \def\cprime{$'$} \def\cprime{$'$} \def\cprime{$'$}
  \def\polhk#1{\setbox0=\hbox{#1}{\ooalign{\hidewidth
  \lower1.5ex\hbox{`}\hidewidth\crcr\unhbox0}}} \def\cdprime{$''$}
  \def\cprime{$'$} \def\cprime{$'$} \def\cprime{$'$} \def\cprime{$'$}
\providecommand{\bysame}{\leavevmode\hbox to3em{\hrulefill}\thinspace}
\providecommand{\MR}{\relax\ifhmode\unskip\space\fi MR }
\providecommand{\MRhref}[2]{%
  \href{http://www.ams.org/mathscinet-getitem?mr=#1}{#2}
}
\providecommand{\href}[2]{#2}

  \bigskip

  Princeton University, Princeton NJ 08544-1000, \

  \email{kollar@math.princeton.edu},
  \email{dpaz@math.princeton.edu}

\end{document}